\def\spr{\mathbb{R}}

\documentclass[reqno, 11pt]{amsart}

\def\Bbb#1{{\fam\msbfam\relax#1}}
\font\fivemsb=msbm5 \font\sevenmsb=msbm7 \font\tenmsb=msbm10
\newfam\msbfam
\textfont\msbfam=\tenmsb \scriptfont\msbfam\sevenmsb
\scriptscriptfont\msbfam\fivemsb

\def\sqbox{{\vcenter{\hrule height .5pt \hbox{\vrule
width .5pt height 5pt \kern 5pt \vrule width .5pt} \hrule height
.5pt}}}
\def\qed{\sqbox}

\usepackage{amsmath}
\usepackage{amsthm}
\usepackage{amssymb}
\newtheorem{theorem}{Theorem}[section]
\newtheorem{corollary}[theorem]{Corollary}
\newtheorem{lemma}[theorem]{Lemma}
\newtheorem{proposition}[theorem]{Proposition}
\newtheorem{observation}[theorem]{Observation}

\newtheorem{definition}[theorem]{Definition}

\begin{document}

\title{Two dimensional Meixner random vectors of class ${\mathcal M}_L$}
\footnotetext[1]{{\em 2000 Mathematics Subject Classifications:}
05E35, 60H40, 46L53.
\newline
\indent \indent {\em Key words and phrases}: commutator,
annihilation operator, creation operator, preservation operator,
Meixner vector of class ${\mathcal M}_L$.}
\author{AUREL I. STAN}
\address{Department of Mathematics \\
The Ohio State University at Marion \\
1465 Mount Vernon Avenue \\
Marion, OH 43302, U.S.A.\\
E-mail: {\em stan.7@osu.edu}}

\maketitle

\begin{abstract}
The paper is divided into two parts. In the first part we lay down
the foundation for defining the joint
annihilation--preservation--creation decomposition of a finite
family of, not necessarily commutative random variables, and show
that this decomposition is essentially unique. In the second part we
show that any two, not necessarily commutative, random variables $X$
and $Y$, for which the vector space spanned by their annihilation,
preservation, and creation operators equipped with the bracket given
by the commutator, forms a Lie algebra, are equivalent, up to an
invertible linear transformation to two independent Meixner random
variables with mixed preservation operators. In particular if $X$
and $Y$ commute, then they are equivalent, up to an invertible
linear transformation to two independent classic Meixner random
variables. To show this we start with a small technical condition
called ``non--degeneracy".
\end{abstract}

\section{Introduction}

This paper continues the work from \cite{sw08}, by considering
examples in dimensions higher than one. The multi--dimensional case
is much more difficult than the one--dimensional one. We will
restrict our attention to the case $d = 2$. The multi--dimensional
case gives us the opportunity to work also in the non--commutative
framework, by considering ``random variables" that do not commute.
In section 2 we introduce the creation, preservation, and
annihilation ($APC$) operators for a finite family of not
necessarily commutative random variables. In section 3, we make the
crucial observation that the vector space spanned by the identity
and $APC$ operators of most of the Meixner random variables,
equipped with the commutator as a bracket, forms a Lie algebra. We
use this observation, to define the notion of ``Meixner random
vector $(X$, $Y)$ of class ${\mathcal M}_L$". We also give two
important examples of such random vectors, one in which $X$ and $Y$
commute, and another in which they do not commute. In section 4,
making a simple and natural assumption, that we call
``non--degeneracy", we classify all Meixner random vectors $(X$,
$Y)$ of class ${\mathcal M}_L$, by showing that through an
invertible linear transformation, they can be reduced to the two
examples from Section 3.

\section{Background}

Let $(H$, $\langle \cdot$, $\cdot \rangle)$ be a Hilbert space over
${\Bbb R}$. Let $X_1$, $X_2$, $\dots$, $X_d$ be $d$ symmetric
densely defined linear operators on $H$. We denote by ${\mathcal A}$
the algebra generated by $X_1$, $X_2$, $\dots$, $X_d$. We assume
that there exists an element $\phi$ of $H$, such that $\phi$ belongs
to the domain of $g$, for any $g \in {\mathcal A}$. We fix such an
element $\phi$ and call it the {\it vacuum vector}.
\begin{definition}
We call any element $g$ of ${\mathcal A}$, a {\em random variable}.
For any $g$ in ${\mathcal A}$, we define:
\begin{eqnarray}
E[g] & := & \langle g\phi, \phi \rangle
\end{eqnarray}
and call the number $E[g]$ the {\em expectation} of the random
variable $g$. Finally, we call the pair $({\mathcal A}$, $\phi)$ a
{\em probability space supported by} $H$.
\end{definition}
The above definition follows the classic GNS representation (see
\cite{s86} and \cite{vdn92}). We restrict ourselves to finitely
generated algebras, but this is not necessary. We always work with
unital algebras ${\mathcal A}$, that means we assume that $I \in
{\mathcal A}$, where $I$ denotes the identity operator of $H$. It is
not hard to see that, for any $g$ in ${\mathcal A}$, there exists a
polynomial $p(x_1$, $x_2$, $\dots$, $x_d)$, of $d$ non--commutative
variables $x_1$, $x_2$, $\dots$, $x_d$, such that $g = p(X_1$,
$X_2$, $\dots$, $X_d)$. We also introduce the following equivalence
relation:
\begin{definition}
Let $({\mathcal A}$, $\phi)$ and $({\mathcal A}'$, $\phi')$ be two
probability spaces supported by two Hilbert spaces $H$ and $H'$, and
let $E$ and $E'$ denote their expectations. Let $X_1$, $X_2$,
$\dots$, $X_d$ be operators from ${\mathcal A}$, and $X_1'$, $X_2'$,
$\dots$, $X_d'$ operators from ${\mathcal A}'$. We say that the
vector random variables $(X_1$, $X_2$, $\dots$, $X_d)$ and $(X_1'$,
$X_2'$, $\dots$, $X_d')$ are {\em moment equal} and denote this fact
by $(X_1$, $X_2$, $\dots$, $X_d) \equiv (X_1'$, $X_2'$, $\dots$,
$X_d')$, if for any polynomial $p(x_1$, $x_2$, $\dots$, $x_d)$ of
$d$ non--commutative variables, we have:
\begin{eqnarray}
E\left[p\left(X_1, X_2, \dots, X_d\right)\right] =
E'\left[p\left(X_1', X_2', \dots, X_d'\right)\right]. \label{equiv}
\end{eqnarray}
\end{definition}
For any non--negative integer $n$, we define the space $F_n$ as
being the set of all vectors of the form $p(X_1$, $X_2$, $\dots$,
$X_d)\phi \in H$, where $p$ is a polynomial of total degree less
than or equal to $n$. It is clear that $F_n$ is a
finite--dimensional subspace of $H$. Being finite dimensional, $F_n$
is a closed subspace of $H$, for all $n \geq 0$. We have:
\begin{eqnarray*}
F_0 \subset F_1 \subset F_2 \subset \cdots \subset H
\end{eqnarray*}
We define $G_0 := F_0$, and for all $n \geq 1$, $G_n := F_n \ominus
F_{n - 1}$, that means $G_n$ is the orthogonal complement of $F_{n -
1}$ into $F_n$. For any $n \geq 0$, we call $G_n$ the {\it
homogenous chaos space of order} $n$ generated by $X_1$, $X_2$,
$\dots$, $X_d$. We also define the space:
\begin{eqnarray*}
{\mathcal H} & := & \oplus_{n = 0}^{\infty}G_n,
\end{eqnarray*}
and call ${\mathcal H}$ the {\it chaos space} generated by $X_1$,
$X_2$, $\dots$, $X_d$. It is not hard to see that ${\mathcal H}$ is
the closure of the space ${\mathcal A}\phi := \{g\phi \mid g \in
{\mathcal A}\}$ in $H$.\\
\par It is also easy to see, based on the fact that $X_1$, $X_2$, $\dots$, $X_d$
are symmetric operators, that we have the following lemma:
\begin{lemma}
For any index $i \in \{1$, $2$, $\dots$, $d\}$ and any non--negative
integer $n$:
\begin{eqnarray}
X_iG_n & \bot & G_k,
\end{eqnarray}
for all $k \neq n - 1$, $n$, and $n + 1$.
\end{lemma}
See also \cite{aks03} for a proof. It follows now, mutatis mutandis
as in \cite{aks03}, that for any $i \in \{1$, $2$, $\dots$, $d\}$,
there exist three operators $a^-(i)$, $a^0(i)$, and $a^+(i)$, called
the {\it annihilation}, {\it preservation}, and {\it creation}
operators, respectively, such that:
\begin{eqnarray}
X_i & = & a^-(i) + a^0(i) + a^+(i). \label{xaaa}
\end{eqnarray}
In (\ref{xaaa}), the domain of $X_i$, $a^-(i)$, $a^0(i)$, and
$a^+(i)$ is understood to be ${\mathcal A}\phi$. It is important to
remember that, for any $i \in \{1$, $2$, $\dots$, $d\}$ and $n \geq
0$, $a^-(i) : G_n \to G_{n - 1}$, $a^0(i) : G_n \to G_n$, and
$a^+(i) : G_n \to G_{n + 1}$, where $G_{-1} := \{0\}$ is the null
space.\\
\par If $Y$ and $Z$ are two operators, then we define their
commutator $[Y$, $Z]$ as:
\begin{eqnarray*}
[Y, Z] & := & YZ - ZY.
\end{eqnarray*}
It is also not hard to see that the operators $X_1$, $X_2$, $\dots$,
$X_d$ commute among themselves if and only if the following
conditions hold, for any $i$, $j \in \{1$, $2$, $\dots$, $d\}$:
\begin{eqnarray}
\left[a^-(i), a^-(j)\right] & = & 0, \label{--}\\
\left[a^-(i), a^0(j)\right] & = & \left[a^-(j), a^0(i)\right], \label{-0}\\
\left[a^-(i), a^+(j)\right] - \left[a^-(j), a^+(i)\right] & = &
\left[a^0(j), a^0(i)\right]. \label{-0+}
\end{eqnarray}
These conditions are derived, as in \cite{aks03} and \cite{aks04},
from the fact that, for all $n \geq 0$:
\begin{eqnarray*}
[X_i, X_j] : G_n \to G_{n - 2} \oplus G_{n - 1} \oplus G_n \oplus
G_{n + 1} \oplus G_{n + 2},
\end{eqnarray*}
by projecting the equality $[X_i, X_j] = 0$ on the spaces $G_{n -
2}$, $G_{n - 1}$, and $G_n$, respectively. Once we know that these
projections are zero, we can conclude that the other two projections
(on $G_{n + 1}$ and $G_{n + 2}$) are also zero, from duality
arguments, since $a^0$ is a symmetric operator,
while $a^+$ is the adjoint of $a^-$.\\
\par If $X_1$, $X_2$, $\dots$, $X_d$ are classic random variables
defined on the same probability space $(\Omega$, ${\mathcal F}$,
$P)$ and having finite moments of any order, then we can take $H =
L^2(\Omega$, ${\mathcal F}$, $P)$ and $\phi = 1$, i.e., the constant
random variable equal to $1$. We then regard $X_1$, $X_2$, $\dots$,
$X_d$ as multiplication operators on the space ${\mathcal A}1
\subset H$, where ${\mathcal A}$ is the algebra of the random
variables of the form $p(X_1$, $X_2$, $\dots$, $X_d)$, where $p$ is
a polynomial of $d$ variables. It is clear
that $X_iX_j = X_jX_i$, for all $1 \leq i < j \leq d$.\\
\par Conditions (\ref{--}), (\ref{-0}), and
(\ref{-0+}) separate the Commutative Probability from the
Non--Commutative one.

\begin{definition} Let $H$ be a Hilbert space,
$({\mathcal A}$, $\phi)$ a probability space supported by $H$, and
$\{X_i\}_{1 \leq i \leq d}$ elements of ${\mathcal A}$, that are
symmetric operators. Let $\{a_{x_i}^-\}_{1 \leq i \leq d}$,
$\{a_{x_i}^0\}_{1 \leq i \leq d}$, and $\{a_{x_i}^+\}_{1 \leq i \leq
d}$ be three families of linear operators, defined on subspaces of
$H$, such that, $\phi$ belongs to the domain of
$a_{i_1}^{\epsilon_1}a_{i_1}^{\epsilon_2} \cdots
a_{i_n}^{\epsilon_n}$, for all $n \geq 1$, $(i_1$, $i_2$, $\dots$,
$i_n) \in \{1$, $2$, $\dots$, $d\}^n$, and $(\epsilon_1$,
$\epsilon_2$, $\dots$, $\epsilon_n) \in \{-$, $0$, $+\}^n$. We say
that these families of operators form {\em a joint
annihilation-preservation-creation} (APC) {\em decomposition of}
$\{X_i\}_{1 \leq i \leq d}$ {\em relative to} ${\mathcal A}$, if the
following conditions hold:
\begin{eqnarray}
X_i & = & a_{x_i}^- + a_{x_i}^0 + a_{x_i}^+ \label{decomposition},\\
\left(a_{x_i}^+\right)^*|{\mathcal A}\phi & = & a_{x_i}^-|{\mathcal
A}\phi
\label{+dual-},\\
a_{x_i}^-H_n & \subset & H_{n - 1}, \label{goesback}\\
a_{x_i}^0H_n & \subset & H_n, \label{0preserves}
\end{eqnarray}
for all $1 \leq i \leq d$ and $n \geq 0$, where $H_{-1} := \{0\}$,
$H_0 := {\Bbb R} \phi$, and $H_k$ is the vector space spanned by all
vectors of the form $a_{x_{i_1}}^+a_{x_{i_2}}^+ \cdots
a_{x_{i_k}}^+\phi$,
where $i_1$, $i_2$, $\dots$, $i_k \in \{1$, $2$, $\dots$, $d\}$, for all $k \geq 1$.\\
We call $a_{x_i}^-$ an {\em annihilation} operator, $a_{x_i}^0$ a
{\em preservation} operator, and $a_{x_i}^+$ a {\em creation}
operator, for all $1 \leq i \leq d.$
\end{definition}
For all $1 \leq i \leq d$, by $\left(a_{x_i}^+\right)^*|{\mathcal
A}\phi = a_{x_i}^-|{\mathcal A}\phi$ we mean:
\begin{eqnarray*}
\langle a_{x_i}^+u, v\rangle & = & \langle u, a_{x_i}^-v\rangle,
\end{eqnarray*}
for all $u$ and $v$ in ${\mathcal A}\phi$, where $\langle \cdot$,
$\cdot \rangle$ denotes the inner product of $H$. Since, for all $i
\in \{1$, $2$, $\dots$, $d\}$, $X_i$ is a symmetric operator, we
conclude from (\ref{+dual-}), that:
\begin{eqnarray}
\left(a_{x_i}^0\right)^*|{\mathcal A}\phi & = & a_{x_i}^0|{\mathcal
A}\phi \label{0dual0}.
\end{eqnarray}

Let us observe that for any $X_1$, $X_2$, $\dots$, $X_d \in
{\mathcal A}$, we can consider the unital algebra ${\mathcal A}'
\subset {\mathcal A}$ generated by $X_1$, $X_2$, $\dots$, $X_d$.
Doing the construction described before, by considering, for all $n
\geq 0$, the space $F_n$ of all vectors of the form $p(X_1$, $X_2$,
$\dots$, $X_d)\phi$, where $p$ is a polynomial of degree at most
$n$, then defining $G_n := F_n \ominus F_{n - 1}$, and so on, we can
construct the annihilation, preservation, and creation operators
$a^-(i)$, $a^0(i)$, and $a^+(i)$ of $X_i$, respectively, where $1
\leq i \leq d$. It is now clear that these operators form a joint
($APC$) decomposition of $\{X_i\}_{1 \leq i \leq d}$ relative to
${\mathcal A}'$. We call this decomposition the {\em minimal joint
(APC) decomposition of} $\{X_i\}_{1 \leq i \leq d}$. We can prove
the following lemma about the uniqueness of the $(APC)$
decomposition, which justifies the choice of the word ``minimal".

\begin{lemma}\label{minimalunique}
Let $\{X_i\}_{1 \leq i \leq d}$ be a family of symmetric random
variables in a non--commutative probability space $({\mathcal A}$,
$\phi)$, and $\{a_{x_i}^-\}_{1 \leq i \leq d}$, $\{a_{x_i}^0\}_{1
\leq i \leq d}$, and $\{a_{x_i}^+\}_{1 \leq i \leq d}$ a joint
($APC$) decomposition of $\{X_i\}_{1 \leq i \leq d}$ relative to
${\mathcal A}$. Let ${\mathcal A}'$ be the algebra generated by
$\{X_i\}_{1 \leq i \leq d}$, and $\{a^-(i)\}_{1 \leq i \leq d}$,
$\{a^0(i)\}_{1 \leq i \leq d}$, and $\{a^+(i)\}_{1 \leq i \leq d}$
the minimal joint ($APC$) decomposition of $\{X_i\}_{1 \leq i \leq
d}$. Then for any $i \in \{1$, $2$, $\dots$, $d\}$ and any $\epsilon
\in \{-$, $0$, $+\}$, we have:
\begin{eqnarray}
a_{x_i}^{\epsilon}|{\mathcal A'}\phi & = & a^{\epsilon}(i)|{\mathcal
A'}\phi. \label{unq}
\end{eqnarray}
Moreover, if ${\mathcal A}''$ denotes the algebra generated by
$\cup_{i = 1}^d\{a_{x_i}^-$, $a_{x_i}^0$, $a_{x_i}^+\}$, then
\begin{eqnarray}
{\mathcal A}''\phi & = & {\mathcal A}'\phi \label{samespace}
\end{eqnarray}
and $\{a_{x_i}^-\}_{1 \leq i \leq d}$, $\{a_{x_i}^0\}_{1 \leq i \leq
d}$, and $\{a_{x_i}^+\}_{1 \leq i \leq d}$ is a joint ($APC$)
decomposition of $\{X_i\}_{1 \leq i \leq n}$ relative to ${\mathcal
A}''$.
\end{lemma}

{\bf Proof.} \ As before, for all $n \geq 0$, let $F_n$ be the space
of all polynomials of $d$ non--commutative variables of degree at
most $n$ (in which the variables $x_1$, $x_2$, $\dots$, $x_n$, are
replaced by $X_1$, $X_2$, $\dots$, $X_n$), $G_n := F_n \ominus F_{n
- 1}$, and $H_n$ the space spanned by all vectors of the form
$a_{x_{i_1}}^+a_{x_{i_2}}^+ \cdots a_{x_{i_n}}^+\phi$, where $i_1$,
$i_2$, $\dots$, $i_n \in \{1$, $2$,
$\dots$, $d\}$. \\
\ \\
\noindent {\bf Claim 1.} \ For all $n \geq 0$, $H_n \subset F_n$.\\
\par We prove this claim by induction on $n$. For $n = 0$, we have
$H_0 = F_0 = {\Bbb R}\phi$. Let us assume that for some $n \geq 1$,
we have $H_k \subset F_k$, for all $k \leq n - 1$, and prove that
$H_n \subset F_n$. We need to show that $a_{x_{i_1}}^+a_{x_{i_2}}^+
\cdots a_{x_{i_n}}^+\phi \in F_n$, for all $i_1$, $i_2$, $\dots$,
$i_n \in \{1$, $2$, $\dots$, $d\}$. Indeed, we have:
\begin{eqnarray*}
& \ & X_{i_1}X_{i_2} \dots X_{i_n}\phi\\
& = & (a_{x_{i_1}}^+ + a_{x_{i_1}}^0 + a_{x_{i_1}}^-)(a_{x_{i_2}}^+
+ a_{x_{i_2}}^0 + a_{x_{i_2}}^-) \cdots (a_{x_{i_n}}^+ +
a_{x_{i_n}}^0 + a_{x_{i_n}}^-)\phi\\
& = & a_{x_{i_1}}^+a_{x_{i_2}}^+ \cdots a_{x_{i_n}}^+\phi +
\sum_{(\epsilon_1, \epsilon_2, \dots, \epsilon_n) \in \{-, 0, +\}^n
\setminus \{(+, + , \cdots,
+)\}}a_{x_{i_1}}^{\epsilon_1}a_{x_{i_2}}^{\epsilon_2} \cdots
a_{x_{i_n}}^{\epsilon_d}\phi.
\end{eqnarray*}
Because at least one operator from the terms of the last sum is a
preservation or annihilation operator, it follows from
(\ref{0preserves}) and (\ref{goesback}), that:
\begin{eqnarray*}
\sum_{(\epsilon_1, \epsilon_2, \dots, \epsilon_n) \in \{-, 0, +\}^n
\setminus \{(+, + , \cdots,
+)\}}a_{x_{i_1}}^{\epsilon_1}a_{x_{i_2}}^{\epsilon_2} \cdots
a_{x_{i_n}}^{\epsilon_d}\phi & \in & H_0 + H_1 + \cdots + H_{n - 1}.
\end{eqnarray*}
Thus, using the induction hypothesis, we get:
\begin{eqnarray*}
& \ & a_{x_{i_1}}^+a_{x_{i_2}}^+ \cdots a_{x_{i_n}}^+\phi\\
& = & X_{i_1}X_{i_2} \dots X_{i_n}\phi - \sum_{(\epsilon_1,
\epsilon_2, \dots, \epsilon_n) \in \{-, 0, +\}^n \setminus \{(+, + ,
\cdots, +)\}}a_{x_{i_1}}^{\epsilon_1}a_{x_{i_2}}^{\epsilon_2} \cdots
a_{x_{i_n}}^{\epsilon_d}\phi\\
& \in & F_n - (H_0 + H_1 + \cdots + H_{n - 1})\\
& \subset & F_n - F_{n - 1}\\
& \subset & F_n.
\end{eqnarray*}

{\bf Claim 2.} \ For all $m \neq n$, $H_m \ \bot \ H_n$ (as
subspaces of $H$).\\
\par We prove by induction on $n$, that for all $n \geq 1$, we have:
$H_n \ \bot \ H_k$, for all $0 \leq k \leq n - 1$.
\par For $n = 1$, we need to prove that $H_1 \ \bot \ H_0$, which
reduces to showing that for all $1 \leq i \leq d$, $a_{x_i}^+\phi \
\bot \ \phi$. If $\langle \cdot $, $\cdot \rangle$ denotes the inner
product of $H$, then it follows from condition (\ref{+dual-}) and
(\ref{goesback}) that:
\begin{eqnarray*}
\langle a_{x_i}^+\phi, \phi\rangle & = & \langle \phi,
a_{x_i}^-\phi\rangle\\
& = & \langle \phi, 0 \rangle\\
& = & 0.
\end{eqnarray*}
Thus $a_{x_i}^+\phi \ \bot \ \phi$.
\par Let us suppose that $H_n \ \bot \ H_k$, for all $0 \leq k \leq n -
1$, and prove that $H_{n + 1} \ \bot \ H_{r}$, for all $0 \leq r
\leq n$. To do this, we need to show that if $r \leq n$, then
\begin{eqnarray*}
a_{x_{i_1}}^+a_{x_{i_2}}^+ \cdots a_{x_{i_{n + 1}}}^+\phi & \bot &
a_{x_{j_1}}^+a_{x_{j_2}}^+ \cdots a_{x_{j_r}}^+\phi,
\end{eqnarray*}
for all $i_1$, $i_2$, $\dots$, $i_{n + 1}$, $j_1$, $j_2$, $\dots$,
$j_r \in \{1$, $2$, $\dots$, $d\}$. Using again the condition
(\ref{+dual-}), because $H_{n + 1}$ and $H_r$ are both contained in
$F_{n + 1}$, and $F_{n + 1} \subset {\mathcal A}\phi$, we get:
\begin{eqnarray*}
& \ & \langle a_{x_{i_1}}^+a_{x_{i_2}}^+ \cdots a_{x_{i_{n +
1}}}^+\phi,
a_{x_{j_1}}^+a_{x_{j_2}}^+ \cdots a_{x_{j_r}}^+\phi\rangle\\
& = & \langle a_{x_{i_2}}^+a_{x_{i_3}}^+ \cdots a_{x_{i_{n +
1}}}^+\phi, a_{x_{i_1}}^-a_{x_{j_1}}^+a_{x_{j_2}}^+ \cdots
a_{x_{j_r}}^+\phi\rangle.
\end{eqnarray*}
From (\ref{goesback}), we know that
$a_{x_{i_1}}^-a_{x_{j_1}}^+a_{x_{j_2}}^+ \cdots a_{x_{j_r}}^+\phi
\in H_{r - 1}$. Since $r \leq n$, we have $r - 1 \leq n - 1$, and
thus, it follows from the induction hypothesis that:
\begin{eqnarray*}
\langle a_{x_{i_2}}^+a_{x_{i_3}}^+ \cdots a_{x_{i_{n + 1}}}^+\phi,
a_{x_{i_1}}^-a_{x_{j_1}}^+a_{x_{j_2}}^+ \cdots
a_{x_{j_r}}^+\phi\rangle & = & 0.
\end{eqnarray*}

\noindent {\bf Claim 3.} \ For all $n \geq 0$, $F_n = H_0 \oplus H_1
\oplus \cdots \oplus H_n$, where ``$\oplus$" denotes the orthogonal sum.\\
\par We know from the previous two claims that $H_0 \oplus H_1
\oplus \cdots \oplus H_n \subset F_n$. On the other hand for any
monomial $X_{i_1}X_{i_2} \dots X_{i_m}\phi$, of degree $m$, where $m
\leq n$, we have:
\begin{eqnarray*}
X_{i_1}X_{i_2} \dots X_{i_m}\phi & = & \sum_{(\epsilon_1,
\epsilon_2, \dots, \epsilon_m) \in \{-, 0, +\}^m}
a_{x_{i_1}}^{\epsilon_1}a_{x_{i_2}}^{\epsilon_2} \cdots
a_{x_{i_n}}^{\epsilon_d}\phi\\
& \in & H_0 \oplus H_1 \oplus \cdots \oplus H_m.
\end{eqnarray*}
Since, these monomials are spanning $F_n$, we get: $F_n \subset H_0
\oplus H_1 \oplus \cdots \oplus H_n$.\\
\par It follows now from Claim 3, that for all $n \geq 0$,
\begin{eqnarray*}
H_n & = & \left[H_0 \oplus H_1 \oplus \cdots \oplus H_n\right]
\ominus \left[H_0 \oplus H_1 \oplus \cdots \oplus H_{n - 1}\right]\\
& = & F_n \ominus F_{n - 1}\\
& = & G_n.
\end{eqnarray*}
Let $n \geq 0$ be fixed and let $w \in H_n$, then for any $1 \leq i
\leq d$, we have:
\begin{eqnarray*}
X_iw & = & a_{x_i}^+w + a_{x_i}^0w + a_{x_i}^-w.
\end{eqnarray*}
Since $a_{x_i}^+w \in H_{n + 1} = G_{n + 1}$, $a_{x_i}^0w \in H_n =
G_n$ (from (\ref{goesback})), and $a_{x_i}^-w \in H_{n - 1} = G_{n -
1}$ (from (\ref{0preserves})), we conclude that $a_{x_i}^+w =
a^+(i)w$, $a_{x_i}^0w = a^0(i)w$, and $a_{x_i}^-w = a^-(i)w$, and
the proof of (\ref{unq}) is complete.\\
\par Moreover,
\begin{eqnarray*}
{\mathcal A}''\phi & = & \cup_{n \geq 0}\left(\oplus_{k =
0}^nH_k\right)\\
& = & \cup_{n \geq 0}F_n\\
& = & {\mathcal A}'\phi.
\end{eqnarray*}
Let us observe now that, for each $1 \leq i \leq d$, $X_i =
a_{x_i}^- + a_{x_i}^0 + a_{x_i}^+ \in {\mathcal A}''$. Thus
$\{X_i\}_{1 \leq i \leq d}$ can be regarded as random variables in
the non--commutative probability space $({\mathcal A}''$, $\phi)$.
It is now clear that $\{a_{x_i}^-\}_{1 \leq i \leq d}$,
$\{a_{x_i}^0\}_{1 \leq i \leq d}$, and $\{a_{x_i}^+\}_{1 \leq i \leq
d}$ is a joint $(APC)$ decomposition of $\{X_i\}_{1 \leq i \leq n}$
relative to ${\mathcal A}''$. \hfill{$\qed$}\\
\par Using the above notations, we close this section by defining the
number operator.
\begin{definition}
We define the {\em number operator} ${\mathcal N}$ as the linear
operator whose domain is ${\mathcal A}'\phi$, defined by the
formula:
\begin{eqnarray}
{\mathcal N}u & = & nu,
\end{eqnarray}
for all $n \geq 0$ and all $u \in H_n.$
\end{definition}

\section{Meixner random vectors of class ${\mathcal M}_L$}
\setcounter{equation}{0}

In the one--dimensional case, a {\em Meixner random variable} $X$,
with infinite support, has the Szeg\"o--Jacobi parameters: $\alpha_n
= \alpha n + \alpha_0$ and $\omega_n = \beta n^2 + (t - \beta)n$,
for all $n \geq 1$, where $\alpha$, $\beta$, and $t$ are real
numbers, such that $\beta > 0$ and $t > 0$. The Meixner random
variables (r.v.), with infinite support, are divided, up to a
re--scaling, into five sub--classes: Gaussian, Poisson, negative
binomial, gamma, and two parameter hyperbolic secant r.v.. A Meixner
random variable $X$, with finite support, has the following
Szeg\"o--Jacobi parameters:
\begin{eqnarray*}
\alpha_n & = & \left\{\begin{array}{ccc} \alpha n + \alpha_0 & {\rm
if} & n \leq k - 1\\
0 & {\rm if} & n \geq k \end{array}\right.
\end{eqnarray*}
and
\begin{eqnarray*}
\omega_n & = & \left\{\begin{array}{ccc} \beta n^2 + (t - \beta)n &
{\rm if} & n \leq k - 1\\
0 & {\rm if} & n \geq k \end{array}\right.,
\end{eqnarray*}
where $k$ is a natural number equal to the number of different
values that $X$ takes on with positive probabilities, and $\alpha$,
$\beta$, and $t$ are real numbers, such that if $k \geq 2$, then $t
> 0$ and $t + \beta(k - 1) > 0$. See \cite{c78}, \cite{m34}, and
\cite{sz75}.\\
\par Thus we can write the class ${\mathcal M}$ of all Meixner
random variables as:
\begin{eqnarray}
{\mathcal M} & := & {\mathcal M}_{u} \cup {\mathcal M}_{b},
\end{eqnarray}
where ${\mathcal M}_{u}$ represents the class of all Meixner random
variables with infinite (unbounded) support, and ${\mathcal M}_{b}$
the class of all Meixner random variables with finite (compact) support.\\
\par Let us consider now a Meixner random variable with infinite
support, having the Szeg\"o--Jacobi parameters: $\alpha_n = \alpha n
+ \alpha_0$ and $\omega_n = \beta n^2 + (t - \beta)n$, for all $n
\geq 1$. This is equivalent to the fact that the commutators of the
terms of its minimal ($APC$) decomposition $X = a^- + a^0 + a^+$
are: $[a^-$, $a^+] = 2\beta {\mathcal N} + tI$, and $[a^-$, $a^0] =
\alpha a^-$, where ${\mathcal N}$ denotes the number operator, i.e.,
the linear operator that maps the $n-$th orthogonal polynomial
$f_n(X)$, generated by $X$, to $nf_n(X)$, for all $n \geq 0$. It
also follows that $a^0 = \alpha{\mathcal N} + \alpha_0I$. Thus, if
$\alpha \neq 0$ or $\alpha = \beta = 0$, then we have $[a^-$, $a^+]
= 2(\beta/\alpha)a^0 + [t - 2\alpha_0(\beta/\alpha)]I$ and $[a^-$,
$a^0] = \alpha a^-$, where $\beta/\alpha := 0$ if $\alpha = \beta =
0$. By duality we also obtain;
\begin{eqnarray*}
[a^0, a^+] & = & \left[(a^0)^*, (a^-)^*\right]\\
& = & [a^-, a^0]^*\\
& = & (\alpha a^-)^*\\
& = & \alpha a^+.
\end{eqnarray*}
This means that the vector space spanned by the operators $a^-$,
$a^0$, $a^+$, and $I$ is closed with respect to taking the
commutator. Thus this vector space equipped with the commutator $[
\cdot $, $ \cdot ]$ forms a Lie algebra. If we define $\beta/\alpha
:= \infty$, for $\alpha = 0$ and $\beta \neq 0$, then we can split
the class ${\mathcal M}_{u}$ into two subclasses: ${\mathcal M}_{u,
f}$ and ${\mathcal M}_{u, \infty}$. ${\mathcal M}_{u, f}$ is the
class of all Meixner random variables, with infinite support, for
which $\alpha \neq 0$ or $\alpha = \beta = 0$. ${\mathcal M}_{u,
\infty}$ is the class of all Meixner random variables, with infinite
support, for which $\alpha = 0$ and $\beta \neq 0$. Thus an element
of ${\mathcal M}_{u, \infty}$ is a random variable having a
symmetric (symmetry about $\alpha_0$) re--scaled two parameter
hyperbolic secant distribution. Moreover, the following proposition
holds:

\begin{proposition}\label{definigmeixner}
The Meixner random variables of class ${\mathcal M}_{u, f}$, are
exactly those random variables $X = a^- + a^0 + a^+$, having finite
moments of all orders and infinite support, for which the vector
space $W$ spanned by $a^-$, $a^0$, $a^+$, and $I$, equipped with the
commutator $[ \cdot $, $ \cdot ]$, forms a Lie algebra, where I
denotes the identity operator. The Gaussian random variables are the
only symmetric (symmetric about a number) random variables, for
which the vector space spanned by $a^-$, $a^+$, and $I$, equipped
with the commutator $[ \cdot $, $ \cdot ]$, forms a Lie algebra.
\end{proposition}

{\bf Proof.} \ ($\Rightarrow$) If $X$ is a Meixner random variable
of class ${\mathcal M}_{u, f}$, then we have already explained why
the vector space $(W $, $[ \cdot $, $ \cdot ])$ forms a Lie algebra.\\
\par ($\Leftarrow$) Let us suppose that $(W$, $[ \cdot $, $ \cdot ])$
is a Lie algebra. Let $\{f_n\}_{n \geq 0}$ be the sequence of
orthogonal polynomials generated by $X$, and let $G_n := {\Bbb R}
f_n$, for all $n \geq 0$. Since $(W $, $[ \cdot $, $ \cdot ])$ is a
Lie algebra, there exist four constants $\alpha$, $c$, $d$, and $e$,
such that:
\begin{eqnarray}
[a^-, a^0] & = & \alpha a^- + ca^0 + da^+ + eI. \label{lie}
\end{eqnarray}
Thus for all $n \geq 0$,
\begin{eqnarray*}
[a^-, a^0]f_n(X) & = & \alpha a^-f_n(X) + ca^0f_n(X) + da^+f_n(X) +
eIf_n(X).
\end{eqnarray*}
This is equivalent to:
\begin{eqnarray*}
(\alpha_n - \alpha_{n - 1})\omega_nf_{n - 1}(X) & = & \alpha
\omega_nf_{n - 1}(X) + (c\alpha_n + e)f_n(X) + df_{n + 1}(X),
\end{eqnarray*}
for all $n \geq 0$, where $\alpha_{-1} := 0$. Since $X$ has an
infinite support, $f_{n - 1}(X)$, $f_n(X)$, and $f_{n + 1}(X)$ are
linearly independent, for all $n \geq 1$. Because $\omega_n > 0$, if
we look at the coefficient of $f_{n - 1}$ in the last equality, we
get $\alpha_n - \alpha_{n - 1} = \alpha$, for all $n \geq 1$. Thus:
\begin{eqnarray*}
\alpha_n & = & \alpha_0 + \sum_{k = 1}^{n}(\alpha_k - \alpha_{k - 1})\\
& = & \alpha_0 + \sum_{k = 1}^{n}\alpha \\
& = & \alpha_0 + \alpha n,
\end{eqnarray*}
for all $n \geq 1$.\\
\par Because $(W$, $[ \cdot $, $ \cdot ])$ is a Lie algebra, there
exist four constants $p$, $q$, $r$, and $s$, such that:
\begin{eqnarray}
[a^-, a^+] & = & pa^- + qa^0 + ra^+ + sI. \label{lie2}
\end{eqnarray}
Since, for all $n \geq 0$, $[a^-$, $a^+]$ and $qa^0 + sI$ map $G_n$
into $G_n$, $pa^-$ maps $G_n$ into $G_{n - 1}$, and $ra^+$ maps
$G_n$ into $G_{n + 1}$, by restricting (\ref{lie2}) to $G_n$, we
get:
\begin{eqnarray}
[a^-, a^+]|G_n & = & (qa^0 + sI)|G_n, \label{rlie2}
\end{eqnarray}
for all $n \geq 0$. If we define $\omega_0 := 0$, then for all $n
\geq 0$,
\begin{eqnarray*}
[a^-, a^+]|G_n & = & (\omega_{n + 1} - \omega_n)I|G_n
\end{eqnarray*}
and
\begin{eqnarray*}
(qa^0 + sI)|G_n & = & (q\alpha_n + s)I|G_n\\
& = & (q\alpha n + q\alpha_0 + s)I|G_n.
\end{eqnarray*}
Hence, it follows from (\ref{rlie2}), that $\omega_{n + 1} -
\omega_n = q\alpha n + q\alpha_0 + s$, for all $n \geq 0$. Since
$\omega_0 = 0$, we obtain:
\begin{eqnarray*}
\omega_n & = & \sum_{k = 0}^{n - 1}(\omega_{k + 1} - \omega_k)\\
& = & \sum_{k = 0}^{n - 1}(q\alpha k + q\alpha_0 + s)\\
& = & q\alpha \frac{n(n - 1)}{2} + (q\alpha_0 + s)n\\
& = & \frac{q\alpha}{2}n^2 + \left(q\alpha_0 + s -
\frac{q\alpha}{2}\right)n,
\end{eqnarray*}
for all $n \geq 1$. Thus $X$ is a Meixner random variable. Moreover,
the numbers $\alpha$ and $\beta$, of $X$, satisfy the relation
$\beta = (q/2)\alpha$, which means that if $\alpha = 0$, then $\beta
= 0$, too. Hence $X$ is of class ${\mathcal M}_{u, f}$.
\hfill{$\qed$}\\

We also have the following proposition:

\begin{proposition}
Let $X$ be a centered random variable, taking on only $k$ different
values with positive probabilities. The space $W$, spanned by the
identity and ($APC$) operators $a_x^-$, $a_x^0$, and $a_x^+$ of $X$,
equipped with the bracket given by the commutator, is a Lie algebra
if and only if $X$ is a Meixner random variable, taking on only $k$
different values, such that, the constants $t$ and $\beta$, involved
in its Szeg\"o--Jacobi parameters, satisfy the following condition:
\begin{eqnarray}
t & = &  -\beta (k - 1), \label{cond}
\end{eqnarray}
and $\alpha \neq 0$. Conditions (\ref{cond}) and $\alpha \neq 0$ are
satisfied exactly by the non--symmetric binomial distributions, with
Krawtchouk orthogonal polynomials.
\end{proposition}

{\bf Proof.} We do the same proof as before, except that due to the
fact that the last non--zero orthogonal polynomial is $f_{k - 1}$,
we have $\omega_k = 0$. Thus $[a^-$, $a^+]f_{k - 1} = -\omega_{k -
1}f_{k - 1}$. Since:
\begin{eqnarray*}
[a^-, a^+]|G_{n - 1} & = & \{2(\beta/\alpha)a^0 + [t -
2\alpha_0(\beta/\alpha)]I\}|G_{n - 1},
\end{eqnarray*}
we must have $[a^-$, $a^+]f_{k - 1} = \{2(\beta/\alpha)a^0 + [t -
2\alpha_0(\beta/\alpha)]I\}f_{k - 1}$. This means
\begin{eqnarray*}
-\omega_{k - 1} & = & 2\beta(k - 1) + t.
\end{eqnarray*}
Since $\omega_{k - 1} = \beta(k - 1)^2 + (t - \beta)(k - 1)$ it
follows
that (\ref{cond}) must hold. \hfill{$\qed$}\\
\par We can split the class of Meixner random variables with
finite support as:
\begin{eqnarray}
{\mathcal M} & = & {\mathcal M}_{b, f} \cup {\mathcal M}_{b,
\infty},
\end{eqnarray}
where ${\mathcal M}_{b, f}$ are the non--symmetric binomials, while
${\mathcal M}_{b, \infty}$ are the symmetric ones.\\
\par Let us now define the following class:
\begin{eqnarray}
{\mathcal M}_{L} & := & {\mathcal M}_{u, f} \cup {\mathcal M}_{b, f}
\end{eqnarray}
and call it the {\em Meixner--Lie class}. We can wrap up our
discussion now in the following lemma:

\begin{lemma}\label{MLlemma}
The real vector space spanned by the identity operator, and the
annihilation, preservation, and creation operators of a random
variable, having finite moments of all orders, equipped with the
bracket given by the commutator, forms a Lie algebra if and only if
that random variable is of Meixner--Lie class.
\end{lemma}

It has been known, see for example \cite{koe96} and \cite{kj98},
that various Meixner classes give rise to Lie algebras like:
${\mathfrak su}(1$, $1)$, ${\mathfrak sl}(2)$, and others.\\

\par We can always assume that $X$ is {\em centered}, that means
$\alpha_0 = E[X] = 0$, since otherwise, we can consider the random
variable $X' = X - E[X]$, whose $\omega$--parameters are the same as
those of $X$, while the $\alpha$--parameters are $\alpha'_n =
\alpha_n - \alpha_0$, for all $n \geq 0$ (in the infinite support case).\\

\par Let us take Lemma \ref{MLlemma}
as a starting point in defining the notion of Meixner probability
distributions on $\mathbb{R}^2$ of class ${\mathcal M}_L$.

\begin{definition}
Let $\mu$ be a probability measure on $\spr^2$ having finite moments
of all orders. Let us denote by $(x$, $y)$ a generic vector in
$\spr^2$. Let $a_x^-$, $a_x^0$, and $a_x^+$ be the annihilation,
preservation, and creation operators generated by the operator $X$
of multiplication by $x$, and $a_y^-$, $a_y^0$, and $a_y^+$ the
annihilation, preservation, and creation operators generated by the
operator $Y$ of multiplication by $y$. We call $\mu$ a {\em Meixner
distribution of class} ${\mathcal M}_L$, if the real vector space
$W$ spanned by the operators $a_x^-$, $a_x^0$, $a_x^+$, $a_y^-$,
$a_y^0$, $a_y^+$, and $I$, equipped with the bracket $[ \cdot $, $
\cdot ]$ given by the commutator, forms a Lie algebra.
\end{definition}

\noindent We extend this definition to the non--commutative case in
the following way:

\begin{definition}
Let $({\mathcal A}$, $\phi)$ be a non--commutative probability space
and $X$ and $Y$ two random variables from ${\mathcal A}$. Let
$\{a_u^-$, $a_u^0$, $a_u^+\}_{u \in \{x, y\}}$ be their minimal
joint ($APC$) decomposition. We say that the pair $(X$, $Y)$ is a
{\em Meixner random vector of class} ${\mathcal M}_L$ if the real
vector space $W$ spanned by the operators $a_x^-$, $a_x^0$, $a_x^+$,
$a_y^-$, $a_y^0$, $a_y^+$, and $I$, equipped with the bracket $[
\cdot $, $ \cdot ]$ given by the commutator, forms a Lie algebra.
\end{definition}

\begin{definition}
Let $X$ and $Y$ be two random variables in a non--commutative
probability space $({\mathcal A}$, $\phi)$ supported by the Hilbert
space $H$. We say that the random vector $(X$, $Y)$ is {\em
non--degenerate} if the vectors $X\phi$, $Y\phi$, and $\phi$ are
linearly independent in $H$. In particular if both $X$ and $Y$ are
centered, then the random vector $(X$, $Y)$ is non--degenerate if
$X\phi$ and $Y\phi$ are linearly independent.
\end{definition}

\begin{proposition}
If $(X$, $Y)$ is a non--degenerate random vector, then the
annihilation operators $a_x^-$ and $a_y^-$, of $X$ and $Y$, are
linearly independent.
\end{proposition}

{\bf Proof.} \ Let $c$ and $d$ be two real numbers such that
\begin{eqnarray*}
ca_x^- + da_y^- & = & 0.
\end{eqnarray*}
Taking the adjoint in both sides of this equality we get:
\begin{eqnarray*}
ca_x^+ + da_y^+ & = & 0.
\end{eqnarray*}
Since $a_x^0 \phi = E[X]\phi$ and $a_y^0\phi = E[Y]\phi$, we have:
\begin{eqnarray*}
& \ & cX\phi + dY\phi - \left(cE[X] + dE[Y]\right)\phi\\
& = & c(a_x^- + a_x^0 + a_x^+)\phi + d(a_y^- + a_y^0 + a_y^+)\phi -
\left(ca_x^0 + da_y^0\right)\phi\\
& = & 0.
\end{eqnarray*}
Thus, because $X\phi$, $Y\phi$, and $\phi$ are linearly independent,
we conclude that $c = d = 0$. Therefore, $a_x^-$ and $a_y^-$ are
linearly independent. \hfill{$\qed$}\\
\par We present now two fundamental examples.\\
\ \\
{\bf Example 1.} \ Let $X$ and $Y$ be two independent centered
Meixner random variables of class ${\mathcal M}_L$ defined on the
same probability space $(\Omega$, ${\mathcal F}$, $P)$. Let $\mu$ be
the joint probability distribution of $X$ and $Y$. We can identify
now $X$ and $Y$ with the multiplication operators on $L^2(\spr^2$,
$\mu)$, by the coordinates $x$ and $y$ of the generic vector $(x$,
$y)$, respectively. In this way $X := a_x^- + a_x^0 + a_x^+$ and $Y
:= a_y^- + a_y^0 + a_y^+$. Since $X$ and $Y$ are independent, we
know from \cite{aks03}, that $[a_x^{\epsilon_1}$, $a_y^{\epsilon_2}]
= 0$, for all $(\epsilon_1$, $\epsilon_2) \in \{-$, $0$, $+\}^2$.
Moreover, one can see that:
\begin{eqnarray}
\left[a_x^-, a_x^+\right] & \in & \spr I + \spr a_x^0,\\
\left[a_y^-, a_y^+\right] & \in & \spr I + \spr a_y^0,\\
\left[a_x^-, a_x^0\right] & \in & \spr a_x^-,\\
\left[a_x^0, a_x^+\right] & \in & \spr a_x^+,\\
\left[a_y^-, a_y^0\right] & \in & \spr a_y^-,\\
\left[a_y^0, a_y^+\right] & \in & \spr a_y^+.
\end{eqnarray}
Hence $(W$, $[ \cdot$, $ \cdot ])$ is a Lie algebra, where $W$ is
the real vector space spanned by $I$, $a_x^-$, $a_x^0$, $a_x^+$,
$a_y^-$, $a_y^0$, and $a_y^+$. Thus $\mu$ is a Meixner distribution
on $\spr^2$ of class ${\mathcal M}_L$, or equivalently $(X$, $Y)$ is
a commutative Meixner random vector of class ${\mathcal M}_L$.\\
\ \\
{\bf Example 2.} \ Let $T$ and $Z$ be two independent centered
Meixner random variables of class ${\mathcal M}_L$ defined on the
same probability space $(\Omega$, ${\mathcal F}$, $P)$, having the
same numbers $\alpha = 1$, $t' = 1$ (we are using $t'$ instead of
$t$, since the letter $t$ we will be used later as a subscript), and
the $\beta$ numbers as follows: $\beta_T := (1/2)(cp + dr)$ and
$\beta_Z := (1/2)(js' + kv)$, where $c$, $p$, $d$, $r$, $j$, $s'$,
$k$, and $v$ are some given real numbers (the choice of these
letters will become more clear in the next section), such that:
\begin{eqnarray*}
cs' + dv & = & 0
\end{eqnarray*}
and
\begin{eqnarray*}
jp + kr & = & 0.
\end{eqnarray*}
In other words, using the usual dot product ``$\cdot$" on $\spr^2$,
and the orthogonal relation ``$\bot$" generated by it, we have:
\begin{eqnarray}
\beta_T & = & \frac{1}{2}(c, d) \cdot (p, r), \label{r1}\\
\beta_Z & = & \frac{1}{2}(j, k) \cdot (s', v), \label{r2}\\
(c, d) & \bot & (s', v), \label{r3}\\
(j, k) & \bot & (p, r) \label{r4}.
\end{eqnarray}

Relations (\ref{r1})--(\ref{r4}) are important, for the following
reason. $\beta_T$ and $\beta_Z$ must be non--negative. We can make
sure that they are not negative by remembering that, the dot product
of two vectors is the product of the length of the vectors times the
cosine of the angle formed by the vectors. Thus we can start with
two non--zero vectors $(c$, $d)$ and $(p$, $r)$ forming an acute
angle. Then we rotate the semi--lines supporting these vectors by
$90^{\circ}$, and choose any two non--zero vectors $(s'$, $v)$ and
$(j$, $k)$ on the rotated semi--lines. In this way we know for sure
that $\beta_T$ and $\beta_Z$ are strictly positive and their sides
are perpendicular,
as required by conditions (\ref{r3}) and (\ref{r4}).\\
\par Let $\mu$ be the joint probability distribution of $T$ and $Z$.
As before, we can identify now $T$ and $Z$ with the multiplication
operators (densely defined) on $H = L^2(\spr^2$, $\mu)$, by the
coordinates $t$ and $z$ of the generic vector $(t$, $z)$. In this
way $T = a_t^- + a_t^0 + a_t^+$ and $Z = a_z^- + a_z^0 + a_z^+$. Let
us consider the following symmetric operators:
\begin{eqnarray*}
X & := & a_t^- + (pa_t^0 + s'a_z^0) + a_t^+
\end{eqnarray*}
and
\begin{eqnarray*}
Y & := & a_z^- + (ra_t^0 + va_z^0) + a_z^+.
\end{eqnarray*}
Let ${\mathcal A}$ be the unital algebra generated by $a_t^-$,
$a_t^0$, $a_t^+$, $a_z^-$, $a_z^0$, and $a_z^+$. Let $\phi := 1$
(the constant random variable, of $H$, equal to $1$). Then
$({\mathcal A}$, $1)$ is a non--commutative probability space
supported by $H$,
and $X$ and $Y$ are random variables in $({\mathcal A}$, $1)$.\\
\ \\
{\bf Claim 1.} \ $(X$, $Y)$ is a Meixner vector of class ${\mathcal
M}_L$.\\
\par Indeed, a joint ($APC$) decomposition of $X$ and $Y$, relative to
${\mathcal A}$ is:
\begin{eqnarray}
a_x^- & = & a_t^-,\\
a_x^0 & = & pa_t^0 + s'a_z^0,\\
a_x^+ & = & a_t^+,\\
a_y^- & = & a_z^-,\\
a_y^0 & = & ra_t^0 + va_z^0,\\
a_y^+ & = & a_z^+.
\end{eqnarray}
Using now Lemma \ref{minimalunique}, we conclude that this ($APC$)
decomposition, restricted to the space ${\mathcal A}'1$, where
${\mathcal A}'$ is the unital algebra generated by $X$ and $Y$, is
the minimal ($APC$) decomposition.\\
\par Since $T$ and $Z$ are linearly independent we have:
\begin{eqnarray*}
\left[a_x^{\pm}, a_y^{\pm}\right] & = & \left[a_t^{\pm},
a_z^{\pm}\right]\\
& = & 0
\end{eqnarray*}
and
\begin{eqnarray*}
\left[a_x^{0}, a_y^{0}\right] & = & 0.
\end{eqnarray*}
Because $T$ is a centered Meixner random variable of class
${\mathcal M}_L$, we have:
\begin{eqnarray*}
\left[a_x^-, a_x^+\right] & = & \left[a_t^-, a_t^+\right]\\
& = & \frac{2\beta_T}{\alpha}a_t^0 + t'I\\
& = & (cp + dr)a_t^0 + I.
\end{eqnarray*}
Multiplying both sides of the relation $a_x^0 = pa_t^0 + s'a_z^0$ by
$c$ and both sides of $a_y^0 = ra_t^0 + va_z^0$ by $d$, and then
adding the two relations together, since $cs' + dv = 0$, we get:
\begin{eqnarray*}
(cp + dr)a_t^0 & = & ca_x^0 + da_y^0.
\end{eqnarray*}
Thus, we get:
\begin{eqnarray*}
\left[a_x^-, a_x^+\right] & = & ca_x^0 + da_y^0 + I \in W.
\end{eqnarray*}
Similarly, we have:
\begin{eqnarray*}
\left[a_y^-, a_y^+\right] & = & ja_x^0 + ka_y^0 + I \in W.
\end{eqnarray*}
We also have:
\begin{eqnarray*}
\left[a_x^-, a_x^0\right] & = & \left[a_t^-, pa_t^0 +
s'a_z^0\right]\\
& = & p\left[a_t^-, a_t^0\right] + s'\left[a_t^-, a_z^0\right]\\
& = & p\alpha a_t^- + s'(0)\\
& = & pa_t^-\\
& = & pa_x^- \in W.
\end{eqnarray*}
Similarly, we have:
\begin{eqnarray*}
\left[a_x^-, a_y^0\right] & = & ra_x^- \in W,\\
\left[a_y^-, a_x^0\right] & = & s'a_y^- \in W,\\
\left[a_y^-, a_y^0\right] & = & va_y^- \in W.
\end{eqnarray*}
By duality (taking the adjoint) we also get that the commutators
between the preservation and creation operators belong to $W$.\\
\par Thus the space $W$, spanned by $I$, $a_x^-$, $a_x^0$, $a_x^+$,
$a_y^-$, $a_y^0$, and $a_y^+$, equipped with the commutator $[
\cdot$, $\cdot ]$ is a Lie algebra. Hence, $(X$, $Y)$ is a Meixner
random vector of class
${\mathcal M}_L$.\\
\ \\
{\bf Claim 2.} \ If $T$ and $Z$ are both not almost surely equal to
$0$, then the centered random vector $(X$, $Y)$ is non--degenerate.\\
\par Indeed, in this case the vacuum vector is $\phi = 1$. Because
$a_t^01 = a_z^01 = 0$ (since $T$ and $Z$ were assumed to be centered
from the beginning), we have $X1 = T1 = T$ and $Y1 = Z1 = Z$. Since
$T$ and $Z$ are independent, they must be linearly independent,
because if $T = \lambda Z$, for some $\lambda \neq 0$, then:
\begin{eqnarray*}
0 & = & E[T]E[Z]\\
& = & E[TZ]\\
& = & E[\lambda Z^2]\\
& = & \lambda E[Z^2]\\
& \neq & 0.
\end{eqnarray*}
We get a contradiction. Thus $X1$ and $Y1$ are linearly independent
and so, the centered random vector $(X$, $Y)$ is non--degenerate.\\
\par Let us observe that $X$ and $Y$ commute if and only if $[a_x^-$,
$a_y^0] = [a_y^-$, $a_x^0]$ ((\ref{-0}) is the only axiom of the
commutative probability that is not guaranteed to hold in this
example), which means $r = s' = 0$, and hence: $X = a_t^- + pa_t^0 +
a_t^+$ and $Y = a_z^- + va_z^0 + a_z^+$ are two independent centered
Meixner random variables, with parameters $\alpha_X = p$, $\alpha_Y
= v$, $\beta_X = \beta_T$, $\beta_Y = \beta_Z$, and $t' = 1$.
Moreover, $X$ is of class ${\mathcal M}_L$, since if $\alpha_X = 0$,
then $\beta_X = (1/2)(cp + dr) = 0$, because $p = \alpha_x = 0$ and
$r = 0$. Similarly $Y$ is of class ${\mathcal M}_L$. Therefore, if
$XY =
YX$, then Example 2 reduces to Example 1.\\
\par Since the annihilation and creation operators of $X$ and $Y$
are the same as the annihilation and creation operators of $T$ and
$Z$, and the preservation operators of $X$ and $Y$ are
superpositions of the preservation operators of $T$ and $Z$, and
because $T$ and $Z$ are independent, we call $X$ and $Y$, from
Example 2, {\em two independent Meixner random variables of class
${\mathcal M}_L$ with mixed preservation operators}.
\begin{observation}
If $(X$, $Y)$ is a non--degenerate Meixner random vector of class
${\mathcal M}_L$, then for any invertible linear transformation $T :
\spr^2 \to \spr^2$, $(X'$, $Y') := T(X$, $Y)$ is also a
non--degenerate Meixner random vector of class ${\mathcal M}_L$. In
particular, if $XY = YX$, then $X'Y' = Y'X'$.
\end{observation}

\section{The classification of non--degenerate centered\\
Meixner random vectors of class ${\mathcal M}_L$}

In this section we show that every non--degenerate centered Meixner
random vector $(X$, $Y)$, of class ${\mathcal M}_L$, can be reduced,
through an invertible linear transformation, to a random vector
equivalent to a vector of two independent centered Meixner random
variables $X'$ and $Y'$, of class ${\mathcal M}_L$, with mixed
preservation operators. In particular if $X$ and $Y$ commute, then
$X'$ and $Y'$ are independent.\\
\par We present first the following lemma.

\begin{lemma}\label{coefficientlemma}
If $(X$, $Y)$ is a Meixner random vector of class ${\mathcal M}_L$,
then the following relations hold:
\begin{enumerate}
\item There exist some real numbers $b$, $c$, $d$, $e$, $f$, $g$, $h$,
$j$, $k$, $p$, $q$, $r$, $s$, $r'$, $s'$, $u$, and $v$, such that:
\begin{eqnarray}
\ \ \ \ \ \ \left[a_x^-, a_x^+\right] & = & bI + ca_x^0 + da_y^0,
\label{-x+x}
\end{eqnarray}
\begin{eqnarray}
\ \ \ \ \ \ \left[a_x^-, a_y^+\right] & = & eI + fa_x^0 + ga_y^0,
\label{-x+y}
\end{eqnarray}
\begin{eqnarray}
\ \ \ \ \ \ \left[a_y^-, a_x^+\right] & = & eI + fa_x^0 + ga_y^0,
\label{-y+x}
\end{eqnarray}
\begin{eqnarray}
\ \ \ \ \ \ \left[a_y^-, a_y^+\right] & = & hI + ja_x^0 + ka_y^0,
\label{-y+y}
\end{eqnarray}
\begin{eqnarray}
\left[a_x^-, a_x^0\right] & = & pa_x^- + qa_y^-, \label{-x0x}
\end{eqnarray}
\begin{eqnarray}
\left[a_x^0, a_x^+\right] & = & pa_x^+ + qa_y^+, \label{0x+x}
\end{eqnarray}
\begin{eqnarray}
\left[a_x^-, a_y^0\right] & = & ra_x^- + sa_y^-, \label{-x0y}
\end{eqnarray}
\begin{eqnarray}
\left[a_y^0, a_x^+\right] & = & ra_x^+ + sa_y^+, \label{0y+x}
\end{eqnarray}
\begin{eqnarray}
\ \left[a_y^-, a_x^0\right] & = & r'a_x^- + s'a_y^-, \label{-y0x}
\end{eqnarray}
\begin{eqnarray}
\ \left[a_x^0, a_y^+\right] & = & r'a_x^+ + s'a_y^+, \label{0x+y}
\end{eqnarray}
\begin{eqnarray}
\left[a_y^-, a_y^0\right] & = & ua_x^- + va_y^-, \label{-y0y}
\end{eqnarray}
\begin{eqnarray}
\left[a_y^0, a_y^+\right] & = & ua_x^+ + va_y^+. \label{0y+y}
\end{eqnarray}

\item
\begin{eqnarray}
\ \left[a_x^-, a_y^-\right] & = & 0, \label{-x-y}
\end{eqnarray}
\begin{eqnarray}
\ \left[a_x^+, a_y^+\right] & = & 0, \label{+x+y}
\end{eqnarray}
\begin{eqnarray}
\left[a_x^0, a_y^0\right] & = & 0. \label{0x0y}
\end{eqnarray}

\item If $XY = YX$, then we can take $r' = r$ and $s' = s$.

\end{enumerate}
\end{lemma}

{\bf Proof.} \ Let $\{G_n\}_{n \geq 0}$ be the homogenous chaos
spaces generated by $X$ and $Y$. \\
\ \\
1. Since, for each $n \geq 0$, the operators $[a_x^-, a_x^+]$,
$[a_x^-, a_y^+]$, and $[a_y^-, a_y^+]$ are mapping $G_n$ into $G_n$,
and because $(W$, $[ \cdot $, $ \cdot ])$ is a Lie algebra, each of
these commutators must be a linear combination of only three
operators from the set $\{I$, $a_x^-$, $a_x^0$, $a_x^+$, $a_y^-$,
$a_y^0$, $a_y^+\}$,  namely, $I$, $a_x^0$, and $a_y^0$. Similarly
since, for each $n \geq 0$, $[a_x^-, a_x^0]$, $[a_x^-, a_y^0]$,
$[a_y^-, a_x^0]$, and $[a_y^-, a_y^0]$ map $G_n$ into $G_{n - 1}$,
each of these commutators must be a linear combination of $a_x^-$
and $a_y^-$. Therefore, there must exist some real numbers $b$, $c$,
$d$, $e$, $f$, $h$, $j$, $k$, $p$, $q$, $r$, $s$, $r'$, $s'$, $u$,
and $v$ such that:
\begin{eqnarray*}
\left[a_x^-, a_x^+\right] & = & bI + ca_x^0 + da_y^0,\\
\left[a_x^-, a_y^+\right] & = & eI + fa_x^0 + ga_y^0,\\
\left[a_y^-, a_y^+\right] & = & hI + ja_x^0 + ka_y^0,\\
\left[a_x^-, a_x^0\right] & = & pa_x^- + qa_y^-,\\
\left[a_x^-, a_y^0\right] & = & ra_x^- + sa_y^-,\\
\left[a_y^-, a_x^0\right] & = & r'a_x^- + s'a_y^-,\\
\left[a_y^-, a_y^0\right] & = & ua_x^- + va_y^-.\\
\end{eqnarray*}

\noindent Taking the adjoint in both sides of (\ref{-x+y}),
(\ref{-x0x}), (\ref{-x0y}), (\ref{-y0x}), and (\ref{-y0y}) we get
(\ref{-y+x}), (\ref{0x+x}), (\ref{0y+x}), (\ref{0x+y}), and
(\ref{0y+y}), respectively.\\
\ \\
2. Since, for all $n \geq 0$, $[a_x^+$, $a_y^+]$ maps $G_n$ into
$G_{n + 2}$, and no operator from the set $\{I$, $a_x^-$, $a_x^0$,
$a_x^+$, $a_y^-$, $a_y^0$, $a_y^+\}$, spanning $W$, has this
property, we have:
\begin{eqnarray*}
[a_x^+, a_y^+] & = & 0.
\end{eqnarray*}
Similarly, we can see that:
\begin{eqnarray*}
[a_x^-, a_y^-] & = & 0.
\end{eqnarray*}
Since, for all $n \geq 0$, $[a_x^0$, $a_y^0]$ maps $G_n$ into $G_n$,
there are three real numbers $\alpha$, $\beta$, and $\gamma$ such
that:
\begin{eqnarray}
[a_x^0, a_y^0] & = & \alpha I + \beta a_x^0 + \gamma a_y^0.
\label{00=0}
\end{eqnarray}
The left--hand side of (\ref{00=0}) is an antisymmetric operator,
due to the fact that:
\begin{eqnarray*}
[a_x^0, a_y^0]^* & = & [a_y^0, a_x^0]\\
& = & -[a_x^0, a_y^0],
\end{eqnarray*}
while the right--hand side is a symmetric operator. Thus:
\begin{eqnarray*}
[a_x^0, a_y^0] & = & 0.
\end{eqnarray*}

\noindent 3. If $XY = YX$, it follows from (\ref{-0}) that $[a_x^-$,
$a_y^0] = [a_y^-$, $a_x^0]$, and thus we can take $r' = r$ and $s' =
s$. \hfill{$\qed$}\\
\par Let us not forget that the random variables $X$ and $Y$ are
assumed to be symmetric operators. We are presenting now two
important lemmas.

\begin{lemma}\label{l1}
If $(X$, $Y)$ is a centered non--degenerate Meixner random vector of
class ${\mathcal M}_L$, then there exists an invertible linear
transformation $T_1 : \spr^2 \to \spr^2$, such that the coefficients
$b'$, $e'$, and $h'$ of the Meixner random vector of class
${\mathcal M}_L$, $(X'$, $Y') := T(X$, $Y)$, from Lemma
\ref{coefficientlemma}, are $b' = h' = 1$ and $|e'| < 1$.
\end{lemma}

{\bf Proof.} \ Since both $X$ and $Y$ are centered we have
$a_x^0\phi = E[X]\phi = 0$ and $a_y^0\phi = E[Y]\phi = 0$. Thus, we
get:
\begin{eqnarray*}
E[X^2] & = & \langle (a_x^+ + a_x^0 + a_x^-)X\phi, \phi\rangle\\
& = & \langle X\phi, a_x^-\phi \rangle + \langle X\phi, a_x^0\phi
\rangle + \langle a_x^-X\phi, \phi \rangle\\
& = & 0 + 0 + \langle a_x^-(a_x^- + a_x^0 + a_x^+)\phi, \phi
\rangle\\
& = & \langle a_x^-a_x^+\phi, \phi \rangle\\
& = & \langle a_x^+a_x^-\phi, \phi \rangle + \langle [a_x^-,
a_x^+]\phi, \phi \rangle\\
& = & 0 + \langle (bI + ca_x^0 + da_y^0)\phi, \phi\rangle\\
& = & b,
\end{eqnarray*}
where $\langle \cdot $, $ \cdot \rangle$ denotes the inner product
of the Hilbert space $H$ supporting the non--commutative probability
space $({\mathcal A}$, $\phi)$, in which $X$ and $Y$ are random
variables. Similarly, we have:
\begin{eqnarray*}
E[Y^2] & = & h
\end{eqnarray*}
and
\begin{eqnarray*}
E[XY] & = & \langle \left[a_x^-, a_y^+\right]1, 1\rangle\\
& = & e.
\end{eqnarray*}
Since $(X$, $Y)$ is non--degenerate we have $X\phi \neq 0$ and
$Y\phi \neq 0$. Because $X$ is a symmetric operator, we have:
\begin{eqnarray*}
b & = & E[X^2]\\
& = & \langle X^2\phi, \phi \rangle\\
& = & \langle X\phi, X\phi \rangle\\
& = & \parallel X\phi \parallel^2\\
& > & 0,
\end{eqnarray*}
where $\parallel \cdot \parallel$ denotes the norm of $H$.
Similarly, we have $h > 0$. Let us consider the invertible linear
change of variables (in fact a simple re--scaling):
\begin{eqnarray*}
(X', Y') & := & \left(\frac{1}{\sqrt{b}}X,
\frac{1}{\sqrt{h}}Y\right).
\end{eqnarray*}
It is easy to see now that the coefficients $b'$ and $h'$ of $(X'$,
$Y')$ are $b' = h' = 1$. Thus by this re--scaling we may assume that
$b = h = 1$. It follows now from the Schwarz' inequality, that:
\begin{eqnarray*}
|e| & = & |E[XY]|\\
& = & |\langle XY\phi, \phi \rangle|\\
& = & |\langle Y\phi, X\phi \rangle|\\
& \leq & \parallel X\phi \parallel \cdot \parallel Y\phi \parallel\\
& = & \sqrt{b} \cdot \sqrt{h}\\
& = & 1.
\end{eqnarray*}
Thus $|e| \leq 1$. We cannot have $|e| = 1$, since if this were
true, then we would have equality is the Schwarz' inequality that we
used, which would imply that $X\phi$ and $Y\phi$ are linearly
dependent, contradicting the fact that $(X$, $Y)$ is
non--degenerate. Hence $|e| < 1$. \hfill{$\qed$}

\begin{lemma}\label{nl1}
Let $(X$, $Y)$ be a non--degenerate centered Meixner random vector
of class ${\mathcal M}_L$. There exists an invertible $2 \times 2$
matrix $T$, with real entries such that for the non--degenerate
centered Meixner random vector $(X_1$, $Y_1) := T(X$, $Y)$ of class
${\mathcal M}_L$, the coefficients $q$, $s$, $r'$, and $u$, from
Lemma \ref{coefficientlemma} are $q = s = r' = u = 0$. Moreover, $b
= h = 1$ and $|e| < 1$.
\end{lemma}

{\bf Proof.} \ From Lemma \ref{l1}, we may assume that: $b = h = 1$
and $|e| < 1$. Let's find first some relations between the values of
the coefficients from Lemma \ref{coefficientlemma}, that hold for
all non--degenerate Meixner random vectors $(X$, $Y)$. \\
\par From the Jacobi identity:
\begin{eqnarray*}
[a_x^-, [a_x^0, a_y^0]] + [a_x^0, [a_y^0, a_x^-]] + [a_y^0, [a_x^-,
a_x^0]] & = & 0,
\end{eqnarray*}
using the fact that $a_x^-$ and $a_y^-$ are linearly independent, we
get:
\begin{eqnarray}
sr' & = & uq \label{sr'=uq}
\end{eqnarray}
and
\begin{eqnarray}
s(s' - p) & = & -q(r - v). \label{a+b=c+d1}
\end{eqnarray}
\par From the Jacobi identity:
\begin{eqnarray*}
[a_y^-, [a_x^0, a_y^0]] + [a_x^0, [a_y^0, a_y^-]] + [a_y^0, [a_y^-,
a_x^0]] & = & 0,
\end{eqnarray*}
we obtain:
\begin{eqnarray}
u(s' - p) & = & -r'(r - v). \label{a+b=c+d2}
\end{eqnarray}
Let $\gamma := s' - p$ and $\delta := r - v$. Relations
(\ref{a+b=c+d1}) and (\ref{a+b=c+d2}) become now:
\begin{eqnarray}
s\gamma = -q\delta \label{delta1}
\end{eqnarray}
and
\begin{eqnarray}
u\gamma = -r'\delta \label{delta2}.
\end{eqnarray}
Formulas (\ref{-x0x}), (\ref{-x0y}), (\ref{-y0x}), and (\ref{-y0y})
can now be written as:

\begin{eqnarray*}
\left[a_x^-, a_x^0\right] & = & pa_x^- + qa_y^-\\
\left[a_x^-, a_y^0\right] & = & (v + \delta)a_x^- + sa_y^-\\
\left[a_y^-, a_x^0\right] & = & r'a_x^- + (p + \gamma)a_y^-\\
\left[a_y^-, a_y^0\right] & = & ua_x^- + va_y^-.
\end{eqnarray*}

\noindent {\bf Case 1.} \ If $q = s = r' = u = 0$, we have nothing
to prove and
are done.\\
\ \\
{\bf Case 2.} \ Let us assume that at least one of the numbers $q$,
$s$, $r'$, and $u$ is not equal to zero.\\
\par We try now to find two random variables $Z_w = \alpha_w X +
\beta_w Y$, $w \in \{1$, $2\}$, where $\alpha_w$ and $\beta_w$ are
real numbers, such that there exist some real constants $\lambda_w$
and $\mu_w$ for which the following relations hold:
\begin{eqnarray}
\left[a_{z_w}^-, a_x^0\right] & = & \lambda a_{z_w}^-, \label{necessary1}\\
\left[a_{z_w}^-, a_y^0\right] & = & \mu a_{z_w}^-,
\label{necessary2}
\end{eqnarray}
and the matrix:
\[
\left(\begin{array}{cc}\lambda_1 & \mu_1\\\lambda_2 &
\mu_2\end{array}\right)
\]
is invertible. We drop for the moment the index $w$. It is easy to
see that if $Z = \alpha X + \beta Y$, then an ($APC$) decomposition
of $Z$ is given by: $a_z^{\epsilon} = \alpha a_x^{\epsilon} + \beta
a_y^{\epsilon}$, for all $\epsilon \in \{-$, $0$, $+\}$. Thus:
\begin{eqnarray*}
\left[a_z^-, a_x^0\right] & = & \alpha\left[a_x^-, a_x^0\right] +
\beta\left[a_y^-, a_x^0\right]\\
& = & (\alpha p + \beta r')a_x^- + [\alpha q + \beta (p +
\gamma)]a_y^-.
\end{eqnarray*}
It follows from here that if we want (\ref{necessary1}) to hold,
then the coefficients $\alpha p + \beta r'$ and $\alpha q + \beta (p
+ \gamma)$ must be proportional to $\alpha$ and $\beta$. Since we
also want $(\alpha$, $\beta) \neq (0$, $0)$, this fact is equivalent
to:
\begin{eqnarray*}
\left|\begin{array}{cc}\alpha p + \beta r' & \alpha q + \beta (p +
\gamma)
\\\alpha & \beta\end{array}\right| & = & 0.
\end{eqnarray*}
This means the following relation must hold:
\begin{eqnarray}
r'\beta^2 - \gamma\beta\alpha - q\alpha^2 & = & 0.
\label{firstequation}
\end{eqnarray}
A similar calculation shows that if we want relation
(\ref{necessary2}) to hold, then we must have:
\begin{eqnarray}
u\beta^2 + \delta\beta\alpha - s\alpha^2 & = & 0.
\label{secondequation}
\end{eqnarray}
Let us look now at the matrix formed by the coefficient vectors
$(r'$, $-\gamma$, $-q)$ and $(u$, $\delta$, $-s)$ of the unknown
vector $(\beta^2$, $\beta \alpha$, $\alpha^2)$ from equations
(\ref{firstequation}) and (\ref{secondequation}). This matrix is:
\[
\left(\begin{array}{ccc} r' & -\gamma & -q\\ u & \delta & -s
\end{array}\right).
\]
Amazingly, it follows from relations (\ref{sr'=uq}), (\ref{delta1}),
and (\ref{delta2}), that the determinant of any $2 \times 2$
sub--matrix of this matrix is zero. Thus the two rows of this matrix
are linearly dependent and so the equations (\ref{firstequation})
and (\ref{secondequation}) are equivalent, unless one of them is the
trivial equation $0 = 0$.\\
\par Since at least one of the numbers $q$, $s$, $r'$, and $u$ is
not zero, we know for sure that at least one of the equations
(\ref{firstequation}) and (\ref{secondequation}) is quadratic in
either $\beta$ or $\alpha$, and the solution(s) of that equation
will also be solution(s) of the other one. Let us assume $r' \neq
0$, and focus on equation (\ref{firstequation}). Let us choose
$\alpha := 1$. Thus this equation becomes:
\begin{eqnarray}
r'\beta^2 - \gamma\beta - q & = & 0. \label{rfirstequation}
\end{eqnarray}
The solutions of this equation are:
\begin{eqnarray}
\beta_1 & = & \frac{\gamma - \sqrt{\gamma^2 + 4r'q}}{2r'}
\end{eqnarray}
and
\begin{eqnarray}
\beta_2 & = & \frac{\gamma + \sqrt{\gamma^2 + 4r'q}}{2r'}.
\end{eqnarray}
We need $\beta_1$ and $\beta_2$ to be real and distinct. If we show
this, then the matrix:
\begin{eqnarray}
\left(\begin{array}{cc} \alpha_1 & \alpha_2\\ \beta_1 & \beta_2
\end{array}\right) & = &
\left(\begin{array}{cc} 1 & 1\\ \beta_1 & \beta_2
\end{array}\right)
\end{eqnarray}
will have a non--zero determinant, and thus the linear
transformation $(X$, $Y) \mapsto (Z_1$, $Z_2)$ will be invertible.\\
\par Thus to achieve our goal we need now to show that the
discriminant:
\begin{eqnarray*}
\Delta & = & \gamma^2 + 4r'q
\end{eqnarray*}
of the quadratic equation (\ref{rfirstequation}) is strictly
positive. To show this we use again one of the Jacobi identities,
applied to the vacuum vector $\phi$, namely:
\begin{eqnarray*}
\left[a_x^-, \left[a_x^0, a_y^+\right]\right]\phi + \left[a_x^0,
\left[a_y^+, a_x^-\right]\right]\phi + \left[a_y^+, \left[a_x^-,
a_x^0\right]\right]\phi & = & 0.
\end{eqnarray*}
Since $a_x^0\phi = a_y^0\phi = 0$, we get:
\begin{eqnarray}
r' - q & = & -e\gamma. \label{req'}
\end{eqnarray}
Let us remember, that from the Schwarz' inequality and the fact that
$(X$, $Y)$ is non--degenerate we know that $|e| < 1$. Thus we have:
\begin{eqnarray}
-4r'q & \leq & (r' - q)^2 \label{square}\\
& = & e^2\gamma^2 \nonumber\\
& \leq & 1\gamma^2. \label{1g}
\end{eqnarray}
It follows now that $\Delta = \gamma^2 + 4r'q \geq 0$. This
inequality must be strict, since if $\Delta = 0$, then we must have
equality in (\ref{square}), which means $(r' + q)^2 = 0$. Thus we
would have $q = -r'$, which would imply: $e\gamma = - 2r' \neq 0$.
Hence $\gamma \neq 0$, and thus (\ref{1g}) is a strict inequality,
which shows that $\Delta > 0$. \\
\par It follows now easily from (\ref{necessary1}) and
(\ref{necessary2}) that, for all $(i$, $j) \in \{1$, $2\}^2$, we
have:
\begin{eqnarray}
\left[a_{z_i}^-, a_{z_j}^0\right] & \in & \spr a_{z_i}^-.
\end{eqnarray}
\par Re--scaling the random variables $Z_1$ and $Z_2$, we may assume
that their coefficients $b$ and $h$ are both equal to $1$.
\hfill{$\qed$}

\begin{lemma}
If $(X$, $Y)$ is a centered non--degenerate Meixner random vector of
class ${\mathcal M}_L$, whose constants from Lemma
\ref{coefficientlemma} satisfy the conditions $b = h = 1$ and $q = s
= r' = u = 0$, then:
\begin{eqnarray}
cs' + dv & = & 0, \label{rel1}\\
jp + kr & = & 0, \label{rel2}\\
fp + gr & = & 0, \label{rel3}\\
fs' + gv & = & 0 \label{rel4}.
\end{eqnarray}
\end{lemma}

{\bf Proof.} \ From the Jacobi identity:
\begin{eqnarray*}
\left[a_x^-, \left[a_y^-, a_x^+\right]\right] + \left[a_y^-,
\left[a_x^+, a_x^-\right]\right] + \left[a_x^+, \left[a_x^-,
a_y^-\right]\right] & = & 0
\end{eqnarray*}
we get:
\begin{eqnarray*}
0 & = & \left[a_x^-, \left[a_y^-, a_x^+\right]\right] +
\left[a_y^-, \left[a_x^+, a_x^-\right]\right]\\
& = & \left[a_x^-, I + fa_x^0 + a_y^0\right] -
\left[a_y^-, I + ca_x^0 + da_y^0\right]\\
& = & (fp + gr)a_x^- -(cs' + dv)a_y^-.
\end{eqnarray*}
Since $(X$, $Y)$ is  non--degenerate, $a_x^-$ and $a_y^-$ are
linearly independent. Thus we obtain that relations (\ref{rel1}) and
(\ref{rel3})hold.\\
\par Similarly, from the Jacobi identity:
\begin{eqnarray*}
\left[a_y^-, \left[a_x^-, a_y^+\right]\right] + \left[a_x^-,
\left[a_y^+, a_y^-\right]\right] + \left[a_y^+, \left[a_y^-,
a_x^-\right]\right] & = & 0,
\end{eqnarray*}
we conclude that (\ref{rel2}) and (\ref{rel4}) hold.
\hfill{$\qed$}.\\

\par We are ready now to present the main theorem.

\begin{theorem}
If $(X$, $Y)$ is a non--degenerate centered Meixner random vector,
then there exists an invertible linear transformation $S : \spr^2
\to \spr^2$, such that the random vector $(X'$, $Y') := S(X$, $Y)$
is equivalent (moment equal) to a random vector of two independent
Meixner random variables with mixed preservation operators of class
${\mathcal M}_L$. Equivalently, the vector space spanned by the
identity operator and the joint ($APC$) operators of $X$ and $Y$ is
isomorphic, as a Lie algebra, to the vector space spanned by the
identity operator and the joint ($APC$) operators of two independent
Meixner random variables of class ${\mathcal M}_L$, with mixed
preservation operators.
\end{theorem}

{\bf Proof.} \ From the previous lemmas we may assume that $q = s =
r' = u = 0$ and $b = h = 1$.\\
\par From the Jacobi identity:
\begin{eqnarray*}
\left[a_y^-, \left[a_x^0, a_x^+\right]\right] + \left[a_x^0,
\left[a_x^+, a_y^-\right]\right] + \left[a_x^+, \left[a_y^-,
a_x^0\right]\right] & = & 0,
\end{eqnarray*}
we get:
\begin{eqnarray*}
\left[a_y^-, pa_x^+\right] + 0 + \left[a_x^+, s'a_y^-\right] & = &
0.
\end{eqnarray*}
This is equivalent to:
\begin{eqnarray}
\gamma \left[a_y^-, a_x^+\right] & = & 0. \label{e1}
\end{eqnarray}
Similarly, from
\begin{eqnarray*}
\left[a_x^-, \left[a_y^0, a_y^+\right]\right] + \left[a_y^0,
\left[a_y^+, a_x^-\right]\right] + \left[a_y^+, \left[a_x^-,
a_y^0\right]\right] & = & 0,
\end{eqnarray*}
we obtain:
\begin{eqnarray}
\delta \left[a_x^-, a_y^+\right] & = & 0. \label{e2}
\end{eqnarray}
Since $\left[a_x^-, a_y^+\right] = \left[a_y^-, a_x^+\right]$, we
conclude from (\ref{e1}) and (\ref{e2}), that if $\gamma \neq 0$ or
$\delta \neq 0$, then:
\begin{eqnarray}
\left[a_x^-, a_y^+\right] & = & 0, \label{zero-x+y}
\end{eqnarray}
and thus we can take $e = f = g = 0$ in (\ref{-x+y}). We analyze now
two cases:\\
\ \\
{\bf Case 1.} \ If $\gamma \neq 0$ or $\delta \neq 0$, then
$[a_x^-$, $a_y^+] = [a_y^-$, $a_x^+] = 0$. We know from the previous
lemma that relations (\ref{rel1}) and (\ref{rel2}) hold. These are
exactly the orthogonality relations from Example 2, from the
previous section. Moreover the commutators of the joint ($APC$)
operators of $X$ and $Y$ are expressed in terms of $I$ and the joint
($APC$) operators of
$X$ and $Y$ in exactly the same way as in Example 1. So we are done.\\
\ \\
\noindent {\bf Case 2.} \ If $\gamma = \delta = 0$, then $p = s'$
and $v = r$. It follows now from the previous lemma that:
\begin{eqnarray}
cp + dv & = & 0, \label{nrel1}\\
jp + dv & = & 0, \label{nrel2}\\
fp + dv & = & 0. \label{nrel3}
\end{eqnarray}
Since $p = s'$, $[a_x^-$, $a_x^0] = pa_x^-$, and $[a_y^-$, $a_x^0] =
s'a_y^-$, we conclude that, for all $z \in \{x$, $y\}$,
\begin{eqnarray*}
\left[a_z^-, a_x^0\right] & = & pa_z^-.
\end{eqnarray*}
By duality it follows now that:
\begin{eqnarray}
\left[a_x^0, a_z^+\right] & = & pa_z^+, \label{0x+z}
\end{eqnarray}
for all $z \in \{x$, $y\}$. Similarly, we have:
\begin{eqnarray}
\left[a_y^0, a_z^+\right] & = & va_z^+, \label{0x+z}
\end{eqnarray}
for all $z \in \{x$, $y\}$.\\
\ \\
{\bf Claim.} \ $a_x^0 = p{\mathcal N}$ and $a_y^0 = v{\mathcal N}$,
where ${\mathcal N}$ denotes the number operator.\\
\par Indeed, for any $n \geq 0$ and any $(z_1$, $z_2$, $\dots$, $z_n)
\in \{x$, $y\}^n$, using the product rule for commutators and the
fact that $a_x^0\phi = 0$, we have:
\begin{eqnarray*}
& \ & a_x^0\left(a_{z_1}^+a_{z_2}^+ \cdots a_{z_n}^+\phi\right)\\
& = & \left(a_{z_1}^+a_{z_2}^+ \cdots a_{z_n}^+\right)a_x^0\phi +
\left[a_x^0, a_{z_1}^+a_{z_2}^+ \cdots a_{z_n}^+\right]\phi\\
& = & \sum_{k = 1}^na_{z_1}^+ \cdots a_{z_{k - 1}}^+ \left[ a_x^0,
a_{z_k}^+\right]a_{z_{k + 1}}^+ \cdots a_{z_n}^+\phi\\
& = & \sum_{k = 1}^na_{z_1}^+ \cdots a_{z_{k - 1}}^+ \left(p
a_{z_k}^+\right)a_{z_{k + 1}}^+ \cdots a_{z_n}^+\phi\\
& = & p\sum_{k = 1}^na_{z_1}^+a_{z_2}^+ \cdots a_{z_n}^+\phi\\
& = & pna_{z_1}^+a_{z_2}^+ \cdots a_{z_n}^+\phi.
\end{eqnarray*}
Hence, for all $n \geq 0$ and all $\xi \in H_n$, we have $a_x^0\xi =
pn\xi$. Since we are interested only on the action of $a_x^0$ on the
space ${\mathcal A}'\phi$, where ${\mathcal A}'$ is the unital
algebra generated by $X$ and $Y$, and ${\mathcal A}'\phi = \cup_{n
\geq 0}\oplus_{k = 0}^nH_k$, we conclude that $a_x^0 = p{\mathcal
N}$. Similarly, we can see that $a_y^0 = v{\mathcal N}$.\\
\ \\
It follows from this claim and the relation (\ref{nrel1}) that:
\begin{eqnarray*}
\left[a_x^-, a_x^0\right] & = & I + ca_x^0 + da_y^0\\
& = & I + (cp + dv){\mathcal N}\\
& = & I + 0{\mathcal N}\\
& = & I.
\end{eqnarray*}
Similarly, it follows from relations (\ref{nrel2}) and (\ref{nrel3})
that:
\begin{eqnarray*}
\left[a_x^-, a_y^0\right] & = & eI
\end{eqnarray*}
and
\begin{eqnarray*}
\left[a_y^-, a_y^0\right] & = & I.
\end{eqnarray*}
Since $|e| < 1$, we can make the following invertible linear change
of variable:
\begin{eqnarray}
\left(X', Y'\right) & = & \left(\frac{1}{\sqrt{2(1 + e)}}(X + Y),
\frac{1}{\sqrt{2(1 - e)}}(X - Y)\right).
\end{eqnarray}
It easy to see now that:
\begin{eqnarray*}
\left[a_{x'}^-, a_{x'}^+\right] & = & I,\\
\left[a_{x'}^-, a_{y'}^+\right] & = & 0,\\
\left[a_{y'}^-, a_{y'}^+\right] & = & I,\\
\left[a_{x'}^-, a_{x'}^0\right] & = & \frac{p + v}{\sqrt{2(1 +
e)}}a_{x'}^-,\\
\left[a_{x'}^-, a_{y'}^0\right] & = & \frac{p - v}{\sqrt{2(1 -
e)}}a_{x'}^-,\\
\left[a_{y'}^-, a_{x'}^0\right] & = & \frac{p + v}{\sqrt{2(1 +
e)}}a_{y'}^-,\\
\left[a_{y'}^-, a_{y'}^0\right] & = & \frac{p - v}{\sqrt{2(1 -
e)}}a_{y'}^-.
\end{eqnarray*}

The proof is now complete for the following reason. Since the Lie
algebra generated by the joint (APC) operators of $(X'$, $Y')$ is
isomorphic to the Lie algebra of the (APC) operators of two
independent Meixner random variables with mixed preservation
operators, and because the joint moments can be recovered from the
commutators, as shown in \cite{sw08}, in the commutative case, and
can easily be extended to the non--commutative case, we conclude
that $(X'$, $Y')$ is moment equal to a random vector whose
components are independent Meixner random variables of class
${\mathcal M}_L$ with mixed preservation operators. \hfill{$\qed$} \
\\
\par Due to the fact that in Example 2, the construction starts with two
independent Meixner random vectors $T$ and $Z$, we can make the
following observation.

\begin{observation}\label{referee3}
The Lie algebra generated by the joint (APC) operators of a two
dimensional non--degenerate Meixner random vector $(X$, $Y)$ of
class ${\mathcal M}_L$ is isomorphic to a subalgebra of
$\mathfrak{sl}(2) \oplus \mathfrak{sl}(2) \oplus \spr I$.
\end{observation}

\begin{corollary}
Let $\mu$ be a Meixner probability measure of class ${\mathcal M}_L$
on $\spr^2$, such that, $\mu$ is not supported by any line of
equation $ax + by = c$, with $a^2 + b^2 > 0$. Then, the following
statements are true:
\begin{enumerate}
\item
Up to an invertible affine transformation, $\mu$ is a product of two
Meixner probability measures of class ${\mathcal M}_L$ on $\spr$.
That means, there exist an invertible linear transformation
${\mathcal S} : \spr^2 \to \spr^2$, a vector $c = (c_1$, $c_2) \in
\spr^2$, and two Meixner probability distributions $\mu_1$ and
$\mu_2$ of class ${\mathcal M}_L$ on $\spr$, such that the measure:
\begin{eqnarray}
\nu(B) & := & \mu({\mathcal S}B + c),
\end{eqnarray}
can be written as:
\begin{eqnarray}
\nu & = & \mu_1 \otimes \mu_2,
\end{eqnarray}
where ${\mathcal S}B + c := \{{\mathcal S}(x$, $y) + c \mid (x$, $y)
\in B\}$, for all Borel subsets $B$ of $\spr^2$.

\item If $\mu$ is not supported by any finite union of
lines of equation $ax + by = c$, with $a^2 + b^2
> 0$, then up to an invertible affine transformation $\mu$ is a
product of two Meixner probability distributions with infinite
support, of class ${\mathcal M}_L$ on $\spr$.
\end{enumerate}
\end{corollary}

\noindent {\bf Final Comment} \ We left out of our discussion the
symmetric two parameter hyperbolic distributions. The Lie Algebra
$W$ that we used in this paper, can be enlarged, so that we can
include also these distributions and characterize the entire Meixner
class using this new Lie Algebra. We hope to do this in another
paper. We would like to mention that the recursive relation among
the orthogonal polynomials, in $d$ variables, that appears in
\cite{p96}, leads to commutation relations among the annihilation,
preservation, and creation operators, that cannot be related to the
present work,
but to the next paper.\\
\ \\
{\bf Acknowledgement} \ The author would like to thank the
anonymous referees and Associate Editor for their kind corrections
and suggestions, which greatly helped him to improve the quality of
this paper. In particular, Observation \ref{referee3} was suggested
by one of the referees, and the author is very grateful for it.

\end{document}